\documentclass[reqno,11pt]{amsart}
\usepackage{amsmath, latexsym, amsfonts, amssymb, amsthm, amscd}
\usepackage{graphics,epsf,psfrag}
\numberwithin{equation}{section}

\newtheorem{theorem}{Theorem}[section]
\newtheorem{lemma}[theorem]{Lemma}
\newtheorem{proposition}[theorem]{Proposition}

\newtheorem{rem}[theorem]{Remark}
\newtheorem{definition}[theorem]{Definition}

\newcommand{\ind}{\mathbf{1}}

\newcommand{\R}{\mathbb{R}}
\newcommand{\Z}{\mathbb{Z}}
\newcommand{\N}{\mathbb{N}}
\renewcommand{\tilde}{\widetilde}

\newcommand{\cL}{{\ensuremath{\mathcal L}} }

\newcommand{\bP}{{\ensuremath{\mathbf P}} }
\newcommand{\bE}{{\ensuremath{\mathbf E}} }

\DeclareMathSymbol{\leqslant}{\mathalpha}{AMSa}{"36} % nicer `smaller or equal'
\DeclareMathSymbol{\geqslant}{\mathalpha}{AMSa}{"3E} % nicer `larger or equal'
\DeclareMathSymbol{\eset}{\mathalpha}{AMSb}{"3F}     % nicer `emptyset'
\newcommand{\dd}{\,\text{\rm d}}             % a straight d for differentials
\newcommand{\maxtwo}[2]{\max_{\substack{#1 \\ #2}}} % max with 2 lines
\newcommand{\sumtwo}[2]{\sum_{\substack{#1 \\ #2}}} % sum with 2 lines
\newcommand{\sumthree}[3]{\sum_{\substack{#1 \\ #2 \\ #3}}} % sum with 3 lines
\newcommand{\limtwo}[2]{\lim_{\substack{#1 \\ #2}}}     % \lim with 2 lines
\newcommand{\bbE}{{\ensuremath{\mathbb E}} }

\newcommand{\bbP}{{\ensuremath{\mathbb P}} }

\newcommand{\go}{\omega}

\makeatletter
\def\captionfont@{\footnotesize}
\def\captionheadfont@{\scshape}

\long\def\@makecaption#1#2{%
  \vspace{2mm}
  \setbox\@tempboxa\vbox{\color@setgroup
    \advance\hsize-6pc\noindent
    \captionfont@\captionheadfont@#1\@xp\@ifnotempty\@xp
        {\@cdr#2\@nil}{.\captionfont@\upshape\enspace#2}%
    \unskip\kern-6pc\par
    \global\setbox\@ne\lastbox\color@endgroup}%
  \ifhbox\@ne % the normal case
    \setbox\@ne\hbox{\unhbox\@ne\unskip\unskip\unpenalty\unkern}%
  \fi
  \ifdim\wd\@tempboxa=\z@ % this means caption will fit on one line
    \setbox\@ne\hbox to\columnwidth{\hss\kern-6pc\box\@ne\hss}%
  \else % tempboxa contained more than one line
    \setbox\@ne\vbox{\unvbox\@tempboxa\parskip\z@skip
        \noindent\unhbox\@ne\advance\hsize-6pc\par}%
\fi
  \ifnum\@tempcnta<64 % if the float IS a figure...
    \addvspace\abovecaptionskip
    \moveright 3pc\box\@ne
  \else % if the float IS NOT a figure...
    \moveright 3pc\box\@ne
    \nobreak
    \vskip\belowcaptionskip
  \fi
\relax
}
\makeatother
\def\writefig#1 #2 #3 {\rlap{\kern #1 truecm
\raise #2 truecm \hbox{#3}}}

\newcommand{\tf}{\textsc{f}}

\newcommand{\bv}{{\underline{v}}}

\begin{document}

\title[Correlation lengths for
pinning models and renewals]{
Correlation lengths for random polymer models and for some
renewal sequences}

 \author{Fabio Lucio Toninelli}
 \address{
 Laboratoire de Physique, UMR-CNRS 5672, ENS Lyon, 46 All\'ee d'Italie,
 69364 Lyon Cedex 07, France
 \hfill\break
 \phantom{br.}{\it Home page:}
 {\tt http://perso.ens-lyon.fr/fabio-lucio.toninelli}}
 \email{fltonine@ens-lyon.fr}
 \date{\today}

\begin{abstract}
  We consider models of directed polymers interacting with a
  one-dimensional defect line on which random charges are placed.
  More abstractly, one starts from renewal sequence on $\Z$ and gives
  a random (site-dependent) reward or penalty to the occurrence of a
  renewal at any given point of $\mathbb Z$.  These models are known
  to undergo a delocalization-localization transition, and the free
  energy $\tf$ vanishes when the critical point is approached from the
  localized region.  We prove that the quenched correlation length
  $\xi$, defined as the inverse of the rate of exponential decay of
  the two-point function, does not diverge faster than $ 1/\tf$.  We
  prove also an exponentially decaying upper bound for the
  disorder-averaged two-point function, with a good control of the
  sub-exponential prefactor.  We discuss how, in the particular case
  where disorder is absent, this result can be seen as a refinement of
  the classical renewal theorem, for a specific class of renewal
  sequences.  
\end{abstract}

\maketitle

\section{Introduction and motivations}

The present work is motivated by the following two problems:

\begin{itemize}

\item {\sl Critical behavior of the correlation lengths for directed
    polymers with (de-)pinning interactions}. Take a homogeneous Markov
    chain $\{S_n\}_{n\ge0}$ on some discrete state space $\Sigma$,
    with $S_0=0$ and law $\bP$. A trajectory of $S$ is interpreted as
    the configuration of a directed polymer in the space $\Sigma\times
    \N$.  In typical examples, $S$ is a simple random walk on
    $\Sigma=\Z^d$ or a simple random walk conditioned to be
    non-negative on $\Sigma=\Z^+$.  Of particular interest is the case
    where the distribution of the first return time of $S$ to zero,
    $K(n):=\bP(\min\{k>0:S_k=0\}=n)$, decays like a power of $n$ for
    $n$ large.  This holds in particular in the case of the simple
    random walks mentioned above.  We want to model the situation
    where the polymer gets a reward (or penalty) $\go_n$ each time it
    touches the line $S\equiv 0$ (which is called {\sl defect
    line}). In other words, we introduce a polymer-line interaction
    energy of the form $$ -\sum_{n=1}^N \go_n \ind_{\{S_n=0\}}, $$
    where $N$ will tend to infinity in the thermodynamic limit. The
    defect line is attractive at points $n$ where $\go_n>0$ and
    repulsive when $\go_n<0$. In particular, one is interested in the
    situation where $\go_n$ are IID quenched random variables.  There
    is a large physics literature (cf. \cite[Chapter 1]{cf:RP} and
    references therein) related to this class of models, due to their
    connection with, e.g., problems of $(1+1)$-dimensional wetting of
    a disordered wall or with the DNA denaturation transition.

In the {\sl localized phase} where the free energy (defined in next
section) is positive and the number of contacts between the polymer
and the defect line, $|\{1\le n\le N: S_n=0\}|$, grows proportionally
to $N$, one knows
\cite{cf:GT_ALEA} that the two-point correlation function
\begin{eqnarray}
\label{eq:2pt}
\left|  \bP_{\infty,\go}(S_{n+k}=0|S_n=0)-  \bP_{\infty,\go}(S_{n+k}=0)
\right|
\end{eqnarray}
decays exponentially in $k$, for almost every disorder
realization. Here, $\bP_{\infty,\go}(.)$ is the Gibbs measure for a
given randomness realization and the index $\infty$ refers to the fact
that the thermodynamic limit has been taken.  
The exponential decay of correlation functions has been 
applied, for instance, to prove sharp results on the maximal excursions
lenght in the localized phase \cite[Theorem 2.5]{cf:GT_ALEA} and
bounds on the finite-size correction to the thermodynamic limit of 
the free energy \cite[Theorem 2.8]{cf:GT_ALEA}.

The inverse of the rate
of decay is identified as a correlation length $\xi$. A natural
question is the relation between $\xi$ and the free energy $\tf$, in
particular in proximity of the delocalization-localization critical
point, where the free energy tends to zero (see next section) and the
correlation length is expected to tend to infinity. The disorder
average of the two-point function \eqref{eq:2pt} is also known 
\cite{cf:GT_ALEA} to decay exponentially with 
$k$, possibly with a {\sl different} rate \cite{cf:T_JSP}.

The important role played by the correlation length, and by its
relation with the free energy, in understanding the critical
properties of disordered pinning models was emphasized in a recent
work by K. Alexander \cite{cf:A}.

\item {\sl Geometric convergence rates for renewal sequences}.
  Consider a renewal sequence $\tau:=\{\tau_i\}_{i=0,1,2,\ldots}$ of
  law $\bP$ defined as follows: $\tau_0=0$, and $\tau_i-\tau_{i-1}$
  are IID  random variables with values in $\N$ and probability
  distribution $p(.)$, where $p(n)\ge 0$ and $\sum_{n\in\N}p(n)=1$.
  The celebrated renewal theorem \cite[Chap. I, Th. 2.2]{cf:asmussen}
  states that
\begin{eqnarray}
  \label{eq:blackw}
  u_n:=\bP(n\in
  \tau)\stackrel{n\to\infty}{\rightarrow}u_\infty:=\frac1{\sum_{n\in
  \N}n p(n)},
\end{eqnarray}
with the convention that $1/\infty=0$. It is natural (and quite useful
in practice, especially in queuing theory applications) to study the
speed of convergence in \eqref{eq:blackw}. In this respect, it is
known (cf. for instance \cite[Chapter VII.2]{cf:asmussen}, \cite{cf:Ney})
that, if
\begin{eqnarray}
  \label{eq:Kend1}
z_{max}:=\sup\{z>0:\sum_{n\in \N}e^{z n}p(n)<\infty\}>0,
\end{eqnarray}
then there exist $r>0$ and $C<\infty$ such that
\begin{eqnarray}
  \label{eq:Kend2}
|u_n-u_\infty|\le C e^{-rn}.
\end{eqnarray}
However, the relation between $z_{max}$ and the largest possible $r$
in Eq.  \eqref{eq:Kend2}, call it $r_{max}$, is not known in
general. A lot of effort has been put in investigating this point, and
in various special cases, where $p(.)$ satisfies some structural
ordering properties, it has been proven that $r_{max}\ge z_{max}$ (see
for instance \cite{cf:BL}, where power series methods are employed and
explicit upper bounds on the prefactor $C$ are given).  In even more
special cases, for instance when $\tau_i$ are the return times of a
Markov chain with some stochastic ordering properties, the optimal
result $r_{max}= z_{max}$ is proved (for details, see
\cite{cf:LT,cf:T_JSP}, which are based on coupling techniques).
However, the equality $r_{max}= z_{max}$ cannot be expected in
general. In particular, if $p(.)$ is a geometric distribution,
$$
p(n)=\frac{e^{-nc}}{e^c-1}
$$
with $c>0$, then one sees that $u_n=u_\infty$
for every $n\in\N$ so that $r_{max}=\infty$, while $z_{max}=c$. On the
other hand, if for instance $p(1)=p(2)=1/2$ and $p(n)=0$ for $n\ge3$,
then $z_{max}=\infty$ while $r_{max}$ is finite.  These and other nice
counter-examples are discussed in \cite{cf:BL}.
\end{itemize}
The two problems are known to be strictly related: indeed, in the
homogeneous situation ($\go_n\equiv const$) the law of the collection
$\{n:S_n=0\}$ of points of polymer-defect contact is given, in the
thermodynamic limit, by a renewal process of the type described above,
with $p(n)$ proportional to $ K(n)e^{-n \tf}$ (cf., for instance,
\cite[Chapter 2]{cf:RP}). In this case, therefore, the free energy
$\tf$ plays the role of $z_{max}$ above.

With respect to the first problem listed above, the main result of
this paper is that, in the limit where $\tf$ tends to zero (i.e., when
the parameters of the model are varied in such a way that the critical
point is approached from the localized phase), the correlation length
$\xi$ is at most of order $1/\tf$, for almost every disorder
realization. An exponentially decaying upper bound, with a good
control of the sub-exponential prefactor, is derived also for the {\sl
disorder average} of the two-point function \eqref{eq:2pt},
cf. Equation \eqref{eq:risultaAV} of Theorem \ref{th:main} and the discussion
in Remark \ref{rimarca}.

As a corollary we obtain the following result for the second problem
above: if the jump law $p(.)$ of the renewal sequence is of the form
$$ p(n)=a_{z_{max}}\frac{L(n)}{n^\alpha}e^{-z_{max}n}, $$ with
$1\le\alpha<\infty$ and $L(.)$ a slowly varying function (not
depending on $z_{max}$), then for $z_{max}$ small one has
that $r_{max}\gtrsim z_{max}$ and $C\lesssim z_{max}^{-c}$ for some
positive constant $c$ (see Theorem
\ref{th:main} and Remarks \ref{rimarca}, \ref{rem:X}
below for the precise statements). In particular, this means that
$|u_n-u_\infty|$ starts decaying exponentially (with rate at least
of order 
$z_{max}$) as soon as $n\gg1/z_{max}$.

\section{Notations and main result}

\label{sec:notations}
We will define our ``directed polymer'' model in an abstract way where
the Markov chain $S$ mentioned in the introduction does not appear
explicitly. In this way the intuitive picture of the Markov chain
trajectory as representing a directed polymer configuration is
somewhat hidden, but the advantage is that the connection with renewal
theory becomes immediate.  The link with the polymer model discussed
in the introduction is made by identifying the renewal sequence $\tau$
below with the set of the return times of the Markov chain $S$ to the
site $0$.

Let $K(.)$ be a probability distribution on $\N:=\{1,2,\ldots\}$, i.e.,
$K(n)\ge0$ for $n\in \N$ and
\begin{eqnarray}
  \label{eq:normalizz}
\sum_{n\in \N}K(n)=1.
\end{eqnarray}
We assume that
  \begin{eqnarray}
    \label{eq:K}
    K(n)=\frac {L(n)}{n^{\alpha}}
  \end{eqnarray}
  for some $1\le \alpha<\infty$. Here, $L(.)$ is a slowly varying
  function, i.e., a positive function $L:\R^+\ni x\rightarrow L(x)\in
  (0,\infty)$ such that $\lim_{x\to\infty}L(xr)/L(x)=1$ for every $r>0$.
Given $x\in\Z$, we construct a renewal process
$\tau:=\{\tau_i\}_{i\in\N\cup \{0\}}$ with law $\bP_x$ as follows:
$\tau_0=x$, and $\tau_i-\tau_{i-1}$ are IID integer-valued random
variables with law $K(.)$.
$\bP_x$ can be naturally seen as a law on the set
$$
\Omega_x:=\{\tau:\tau\subset (\Z\cap [x,\infty))\;\;\mbox{and}\;\; x\in \tau\}.
$$
Note that, thanks to \eqref{eq:normalizz},
$\tau$ is a {\sl recurrent} renewal process (possibly, null-recurrent).

Now we modify the law of the renewal by switching on a random
interaction as follows.  We let $\{\go_n\}_{n\in\Z}$ be a sequence of
IID centered random variables with law $\bbP$ and $\bbE\, \go_0^2=1$.
For simplicity, we require also $\go_n$ to be bounded.  Then, given
$h\in \R$, $\beta\ge0$, $x,y\in \Z$ with $x<y$ and a realization of
$\go$ we let
\begin{eqnarray}
  \label{eq:Boltz}
  \frac{\dd\bP_{x,y,\go}}{\dd\bP_{x}}(\tau)=\frac{e^{\sum_{n=x+1}^y(\beta
\go_n-h)\ind_{\{n\in\tau\}}} } {Z_{x,y,\go}}\ind_{\{y\in \tau\}}
\end{eqnarray}
where, of course,
\begin{eqnarray}
  \label{eq:Z}
  Z_{x,y,\go}=\bE_x\left(e^{\sum_{n=x+1}^y(\beta
  \go_n-h)\ind_{\{n\in\tau\}}}\ind_{\{y\in \tau\}}\right)
\end{eqnarray}
and $\bP_{x,y,\go}$ is still a law on $\Omega_x$.  Note that the
normalization condition \eqref{eq:normalizz} is by no means a
restriction: if we had $\Sigma:=\sum_{n\in \N}K(n)<1$, we could
perform the replacements $K(.)\to K(.)/\Sigma$, $h\to h-\log \Sigma$
in \eqref{eq:Boltz} and the measure $\bP_{x,y,\go}$ would be
unchanged.

One defines the free energy as
\begin{eqnarray}
  \label{eq:F}
  \tf(\beta,h)=\lim_{N\to\infty} \frac 1{2N} \log Z_{-N,N,\go}.
\end{eqnarray}
The convergence holds almost surely and in $L^1(\bbP)$, and
$\tf(\beta,h)$ is $\bbP(\dd \go)$-a.s. constant (see
\cite[Chap. 4]{cf:RP} and \cite{cf:AS}).  It is known that
$\tf(\beta,h)\ge0$: to realize this, it is sufficient to observe that
\begin{eqnarray}
  \label{eq:F>0}
 && \frac 1{2N} \log Z_{-N,N,\go}\ge \frac 1{2N} \log \bE_{-N}
\left(e^{\sum_{n=-N+1}^N(\beta
\go_n-h)\ind_{\{n\in\tau\}}}\ind_{\{\tau_1=N\}}\right)\\ &&=
\frac{\beta \go_N-h}{2N}+\frac1{2N}\log K(2N)
\end{eqnarray}
which tends to zero for $N\to\infty$.  One then
decomposes the phase diagram into {\sl localized} and {\sl
delocalized} regions defined as
\begin{eqnarray}
  \label{eq:loc}
  \mathcal L:=\{(\beta,h):\tf(\beta,h)>0\}
\end{eqnarray}
\begin{eqnarray}
  \label{eq:deloc}
  \mathcal D:=\{(\beta,h):\tf(\beta,h)=0\},
\end{eqnarray}
separated by the critical line
\begin{eqnarray}
  \label{eq:hcrit}
  h_c(\beta):=\inf\{h:\tf(\beta,h)=0\}.
\end{eqnarray}
By convexity, the free energy is continuous in $\beta$ and $h$ and therefore
tends to zero when the critical line is approached from the localized region.
It is known that  typical configurations $\tau$ are very different
in the two regions. Roughly speaking, if $(\beta,h)\in \mathcal L$
then $\tau$ has a finite density of points in $\N$, i.e., for $N$
large
\begin{eqnarray}
\frac1N{\left|\tau\cap\{1,\ldots,N\}\right|}\sim   -\partial_h\tf(\beta,h)>0.
\end{eqnarray}
On the other hand, in $\mathcal D$ the density tends to zero with $N$:
\begin{eqnarray}
\frac1N{\left|\tau\cap\{1,\ldots,N\}\right|}
\left\{
\begin{array}{lcr}
\le (\log N)/N &\mbox{if} &  h>h_c(\beta)\\
\le N^{-1/3}\log N &\mbox{if} & h=h_c(\beta)
  \end{array}
\right.
\end{eqnarray}
(for precise statements see, respectively, \cite[Theorem 1.4, part
(2)]{cf:GT_PTRF} and \cite[Theorem 3.1]{cf:T_JSP}).

Another quantity which will play an important role in the following is
\begin{eqnarray}
  \label{eq:mu}
  \mu(\beta,h)=-\lim_{N\to\infty} \frac 1{2N} \log \bbE \frac1{Z_{-N,N,\go}}.
\end{eqnarray}
As it is known (cf. \cite[Theorem 2.5 and Appendix B]{cf:GT_ALEA})
for $(\beta,h)\in \mathcal L$ one has
\begin{eqnarray}
  \label{eq:ineqmu}
  0<\mu(\beta,h)<\tf(\beta,h),
\end{eqnarray}
while $\tf(\beta,h)=\mu( \beta,h)=0$ in $\mathcal D$. On the other
hand, it is unknown whether the ratio $\tf(\beta,h)/\mu(\beta,h)$
remains bounded for $h\to h_c(\beta)$.  $\mu(\beta,h)$ is related to
the maximal excursion length in the localized phase, $$
\Delta_N:=\maxtwo{0<i<j<N:}{\{i,\ldots,j\}\cap
\tau=\emptyset}\left|j-i\right|, $$ in the sense that essentially
$\Delta_N\simeq \log N/\mu(\beta,h)$, see \cite[Theorem
2.5]{cf:GT_ALEA} (cf. also \cite{cf:AZ} for a proof of the same fact
in a related model, the {\sl heteropolymer at a selective interface}).

As was proven in
\cite{cf:GT_ALEA} (but see also \cite{cf:BidH} for the proof of the almost
sure existence of the
infinite-volume Gibbs measure for the heteropolymer model in the localized
phase), the limit
\begin{eqnarray}
  \label{eq:inftvol}
\bE_{\infty,\go}(f):=  \limtwo{x\to-\infty}{y\to\infty}\bE_{x,y,\go}(f)
\end{eqnarray}
exists, $\bbP(\dd\go)-$a.s., for every $(\beta,h)\in \mathcal L$ and
for every bounded local observable $f$, and is independent of the way
the limits $x\to-\infty$, $y\to\infty$ are performed.  A bounded local
observable is a bounded function $f:\{\tau:\tau\subset \Z\}\rightarrow
\R$ for which there exists $I$,  finite subset of $\Z$, such that
$$ f(\tau_1)=f(\tau_2) $$ whenever $\tau_1\cap I=\tau_2\cap I$.  The
smallest possible $I$ is called support of $f$.  An example of local
observable is $|\{\tau\cap I\}|$, the number of points of $\tau$ which
belong to $I$. On the other hand, $\tau_1$ is not a local observable.

A useful identity is the following: let $a\in\Z$ and $f,g$ be two
local observables, whose supports are contained in
$\{\ldots,a-2,a-1\}$ and $\{a+1,a+2,\ldots\}$, respectively. Then, if
$x<a<y$,
\begin{eqnarray}
  \label{eq:identita}
  \bE_{x,y,\go}(f\,g|a\in\tau)=\bE_{x,a,\go}(f)\bE_{a,y,\go}(g).
\end{eqnarray}
In other words, conditioning on the event that $a$ belongs to $\tau$
makes the process to the left and to the right of $a$
independent. This is easily checked from the definition
\eqref{eq:Boltz} of the Boltzmann-Gibbs measure and from the IID
character of $\tau_i-\tau_{i-1}$ under $\bP_x$.

Our first  result is an exponentially decaying upper bound on the
disorder-averaged  two-point correlation function,
in the localized phase:
\begin{theorem}
\label{th:main}
 Let $\epsilon>0$ and $(\beta,h)\in \mathcal L$. There exists
$C_1:=C_1(\epsilon,\beta,h)>0$ such that, for every $k\in\N$,
\begin{eqnarray}
  \label{eq:risultaAV}
\bbE \left|\bP_{\infty,\go}(k\in
\tau|0\in\tau)- \bP_{\infty,\go}(k\in \tau)\right|\le
\frac1{C_1\mu(\beta,h)^{1/C_1}}
\exp\left(-k \,C_1\,\mu(\beta,h)^{1+\epsilon}\right).
\end{eqnarray}
The constant $C_1(\epsilon,\beta,h)$ does not vanish at the
critical line: for every bounded subset $B\subset \cL$ one has
$\inf_{(\beta,h)\in B}C_1(\epsilon,\beta,h)\ge C_1(B,\epsilon)>0$.
\end{theorem}
\begin{rem}\rm
\label{rimarca}
 Note that Theorem \ref{th:main} is more than just a bound on the rate
 of exponential decay of the disorder-averaged two-point
 correlation. Indeed, thanks to the explicit bound on the prefactor in
 front of the exponential, Eq. \eqref{eq:risultaAV} says that the
 exponential decay, with rate at least of order $\mu^{1+\epsilon}$,
 commences as soon as $k\gg \mu^{-1-\epsilon}|\log \mu|$. This
 observation reinforces the meaning of Eq. \eqref{eq:risultaAV} as an
 upper bound on the correlation length of disorder-averaged
 correlations functions.
\end{rem}
It would be possible, via the Borel-Cantelli Lemma, to extract from
Eq. \eqref{eq:risultaAV} the almost-sure exponential decay of the
disorder-dependent two-point function.  However, from \cite{cf:T_JSP}
one expects the almost-sure exponential decay to be related to $\tf(\beta,h)$
rather than to $\mu(\beta,h)$. Indeed, we have the following:
\begin{theorem}
 \label{th:mainAS}
 Let $\epsilon>0$ and $(\beta,h)\in \mathcal L$. One has for every $k\in\N$
\begin{eqnarray}
\label{eq:risultaAS}
 \left|\bP_{\infty,\go}(k\in
 \tau|0\in\tau)- \bP_{\infty,\go}(k\in \tau)\right|\le
C_2(\go)
\exp\left(- k\,C_1\tf(\beta,h)^{1+\epsilon}\right),
\end{eqnarray}
where $C_1$ is as in Theorem \ref{th:main}, while $C_2(\go):=
C_2(\go,\epsilon,\beta,h)$ is an almost surely finite random variable.
\end{theorem}
Recalling that $\tf>\mu$, it is clear that Theorem \ref{th:mainAS}
cannot be deduced from Theorem \ref{th:main}.
\begin{rem}\rm
It is quite tempting to expect that, in analogy with Theorem
\ref{th:main}, the (random) prefactor $C_2(\go)$ is bounded above by
$$
\frac{C_5(\go,\epsilon,\beta,h)}{\tf(\beta,h)^{C_5(\go,\epsilon,\beta,h)}},
$$ for some random variable $C_5$ such that, say, $\bbE
C_5(\go,\epsilon,\beta,h)\le c(B,\epsilon)<\infty$ for $(\beta,h)$
belonging to a bounded set $B\subset \cL$. This would mean that the
almost sure exponential decay with decay rate at least
of order $\tf^{1+\epsilon}$ commences as soon as $k\gg
n(\go)\tf^{-1-\epsilon}|\log \tf|$, with $n(\go)$ random but
typically of order one even close to the critical point. However, this
kind of result seems to be out of reach with the present techniques.
\end{rem}

\begin{rem}
  \rm
\label{rem:extract}
As can be extracted from the proof of Theorems \ref{th:main} and
\ref{th:mainAS} (see in particular Remark \ref{rem:XY}), if the
slowly varying function $L(n)$ in \eqref{eq:K} tends to a constant for
$n\to\infty$, then one can replace the r.h.s. of Eqs.
\eqref{eq:risultaAV}, \eqref{eq:risultaAS} by
$$
\frac1{C_3(\beta,h)\mu(\beta,h)^{1/C_3(\beta,h)}}
\exp\left(
-k \,C_3(\beta,h)\,\frac{\mu(\beta,h)}{|\log \mu(\beta,h)|}
\right)
$$
and
$$
C_4(\go,\beta,h)
\exp\left(- k\,C_3(\beta,h)\frac{\tf(\beta,h)}{|\log \tf(\beta,h)|}\right)
$$
respectively, with
$\inf_{(\beta,h)\in B}C_3(\beta,h)\ge C_3(B)>0$ and $C_4$ almost surely finite.
\end{rem}

Once the exponential decay of the two-point function is proven, it is
not difficult to obtain similar results for the correlation between any
two given local observables (cf. Remark \ref{rem:decadAB} below for
some more details):
\begin{theorem}
\label{th:AB}
Let $A$ and $B$ be two bounded local observables, with supports $S_A$
and $S_B$, respectively. Assume that $S_A$ is contained in $\Z\cap
(-\infty,0]$ and $S_B\subset \Z\cap [k,\infty)$.  Let
$(\beta,h)\in\cL$, while $\epsilon>0$.  Then,
\begin{eqnarray}
\bbE\,\left|\bE_{\infty,\go}(A B)-\bE_{\infty,\go}(A )\bE_{\infty,\go}(B)
\right|\le \frac{||A||_\infty ||B||_\infty}{C_1 \mu(\beta,h)^{1/C_1}}
\exp\left(-k \,C_1\,\mu(\beta,h)^{1+\epsilon}\right)
\end{eqnarray}
and
\begin{eqnarray}
\left|\bE_{\infty,\go}(A B)-\bE_{\infty,\go}(A )\bE_{\infty,\go}(B)
\right|\le ||A||_\infty ||B||_\infty C_2(\go)
\exp\left(-k \,C_1\,\tf(\beta,h)^{1+\epsilon}\right),
\end{eqnarray}
where $C_1$ and $C_2$ are as in Theorems \ref{th:main} and \ref{th:mainAS}.
\end{theorem}

\section{Sketch of the idea: auxiliary Markov process and coupling}

In this section, we give an informal sketch of the basic ideas
underlying the proof of the upper bounds for the two-point
function. The actual proof is somewhat involved and takes Sections
\ref{sec:markov} to \ref{sec:largeT}.

The basic trick is to associate to the renewal probability $K(.)$ a
Markov process $\{S_t\}_{t\ge x}$ such that, very roughly speaking,
its trajectories are continuous ``most of the time'' and the random
set of times $\{t\in\Z\cap [x,\infty):S_t=0\}$ has the same
distribution as the {\sl discrete} renewal process $\{\tau_i\}_{i\in
\N\cup\{0\}}$ associated to $K(.)$, with law $\bP_{x}$. This
construction is done in Section \ref{sec:markov}, where we see that
$S_.$ is strictly related to the Bessel process \cite{cf:RY} of
dimension $2(\alpha+1)$.  Once we have $S_.$, we switch on the
interaction
$$
-\sum_{n=x+1}^y (\beta\go_n-h)\ind_{\{S_n=0\}}
$$ and in the thermodynamic limit $x\to-\infty,y\to\infty$ we obtain a
new measure $\hat \bP_{\infty,\go}$ on the paths
$\{S_t\}_{t\in\R}$. An important point will be that the process $S_.$,
under $\hat \bP_{\infty,\go}$, is still Markovian, and that the
marginal distribution of $\tau:=\{t\in\Z:S_t=0\}$ is just the measure
$\bP_{\infty,\go}$ defined in Eq.  \eqref{eq:inftvol}.  At that point,
we take two copies $(S^1_., S^2_.)$ of the process, distributed
according to the product measure $\hat \bP^{\otimes 2}_{\infty,\go}$,
and we define the coupling time $\mathcal T(S^1,S^2)=\inf\{t\ge 0:
S^1_t=S^2_t\}$.  From the Markov property it follows that
\begin{eqnarray}
\label{eq:scheccio}
  \left|\bP_{\infty,\go}(k\in \tau|0\in\tau)- \bP_{\infty,\go}(k\in
\tau)\right|\le \hat \bP^{\otimes 2}_{\infty,\go}(\mathcal T(S^1,S^2)>
k|S^1_0=0).
\end{eqnarray}
Indeed, if the two paths meet before time $k$, we can let them proceed
together from then on and they will either both touch zero at $t=k$,
or both will not touch it.  Note that at the left-hand side of
\eqref{eq:scheccio} we have just the quantity we wish to bound in
Theorems \ref{th:main} and \ref{th:mainAS}.  Finally, in order to
prove Eq.  \eqref{eq:risultaAS}, we will show in Section
\ref{sec:ritorni} that, roughly speaking, in the time interval $[0,k]$
two typical (with respect to $\hat \bP^{\otimes 2}_{\infty,\go}$)
configurations of the paths $S^1_.,S^2_.$ come close to each other at
least approximately $ k\, \tf(\beta,h)$ times.  The inequality
\eqref{eq:risultaAS} then follows by estimating what is the
probability that the two (independent!)  paths actually succeed in
avoiding each other every time they are close: it is rather intuitive
that this probability should decrease with $k$ like $\exp(-k
\tf(\beta,h))$.  This explains result \eqref{eq:risultaAS} (forget for
the moment about $\epsilon$).  Inequality \eqref{eq:risultaAV} is
somewhat less intuitive and we do not try to give a heuristic
justification here.  The technical difficulties one meets in turning
this heuristics into a proof are reflected in the necessity of taking
$\epsilon>0$ in Theorem \ref{th:main}.

%\section{Discussion and open problems}

The most natural question left open by our result is whether lower
bounds on the two-point correlation function, complementary to the
upper bounds of Eqs. \eqref{eq:risultaAS}, \eqref{eq:risultaAV}
hold.
In Ref. \cite{cf:T_JSP} a sharp result was proven in a specific case:
if $\bP$ is the law of the zeros of the one-dimensional simple random
walk conditioned to be non-negative (but that proof works also for the
unconditioned simple random walk), then the limit in
\eqref{eq:risultaAS} exists for $(\beta,h)\in\mathcal L$ and equal
exactly $\tf(\beta,h)$.  Similarly, for the disorder-averaged
two-point function the analogous limit exists and equals
$\mu(\beta,h)$. The simplification that occurs in the situation
considered in \cite{cf:T_JSP} is that two trajectories of the Markov
chain which is naturally associated to $K(.)$, i.e., of the simple
random walk, must necessary meet whenever they cross each other. This
avoids the construction of the auxiliary Markov chain and makes the
coupling argument much more efficient.

Let us emphasize that, in general, it is not even proven that the rate
of exponential decay of the (averaged or not) two-point correlation
function tends to zero when the critical point is approached (although
this is very intuitive, and known for instance in the case considered
in \cite{cf:T_JSP}, as  already mentioned).

\section{The Markov process}

\label{sec:markov}

For $\delta\in(2,\infty)$ let $\{\rho^{(s)}_t\}_{ t\ge s}$ be the
Bessel process of dimension $\delta$ and denote its law by
$P_\rho^{(s)}$. The Bessel process is actually well defined also
for $\delta\le2$, but we will not need that here.
For the application we have in mind, we choose the
initial condition $\rho^{(s)}_s=1$. For general properties of the
Bessel process, we refer to \cite[Sections VI.3 and XI.1]{cf:RY}.
This is a diffusion on $\R^+$ with infinitesimal generator
\begin{eqnarray}
  \label{eq:generator}
  \frac12 \frac{d^2}{dx^2}+\frac{\delta-1}{2x}\frac d{dx}.
\end{eqnarray}
For every real $\delta>2$, $\rho^{(s)}_.$ is a transient Markov
process with continuous trajectories (and, if $\rho_s^{(s)}=0$ were
chosen as initial condition, for  $\delta$ integer $\rho^{(s)}_t$ would
have the same law as the absolute value of the standard Brownian
motion in $\R^\delta$ started at the origin at time $s$). The
transition semi-group associated to $\rho^{(s)}_.$, which gives the
probability of being in $y$ at time $t$ having started at $x$ at time
$0$, is known explicitly \cite{cf:RY}: its density in $y$ with respect
to the Lebesgue measure is given, for $t,x>0$, by
\begin{eqnarray}
  \label{eq:semigroup}
  p_t^\delta(x,y):=\frac yt\left(\frac yx\right)^\nu
e^{-(x^2+y^2)/(2t)}I_\nu\left(\frac{xy}t\right)
\end{eqnarray}
where $\nu:=(\delta/2)-1$ and $I_.(.)$ is the modified Bessel function of first
kind \cite[Chapter 7.2.2]{cf:erdely}.

Recall our choice $\rho^{(s)}_s=1$ and define
$T^{(s)}:=\inf\{t> s:\rho^{(s)}_t=1/2\}$.
(As will be clear from the proof, the values $1$ and $1/2$ could be
replaced by any $a,b$ with $a>b>0$.)  Then,
$0<P_\rho^{(s)}(T^{(s)}<\infty)<1$, the upper bound being a
consequence of transience. We let also $\{\hat \rho^{(s)}_t\}_{t\ge
s}$ with law $\hat P_\rho^{(s)}$ be the process $\rho^{(s)}_.$
conditioned on $T^{(s)}<\infty$. Finally, for $n\in\N$ we set
$\mathcal K^{(\delta)}(n):=\hat P_\rho^{(0)}(T^{(0)}\in(n-1,n])$ so
that
\begin{eqnarray}
  \label{eq:Kbess1}
\sum_{n\in\N}\mathcal K^{(\delta)}(n)=1.
\end{eqnarray}
One can prove (cf. Appendix \ref{sec:appBessel}; the proof is an
immediate consequence of results in \cite{cf:kent} and
\cite{cf:mourad}) that
\begin{eqnarray}
  \label{eq:Kbessel}
\lim_{n\to\infty}n^{\delta/2} \mathcal  K^{(\delta)}(n) \in (0,\infty),
\end{eqnarray}
the existence of the limit being part of the statement.

Note that $\hat \rho^{(s)}_.$ {\sl is not} a Markov process.
Indeed, for instance,
\begin{eqnarray}
&&\hat P_\rho^{(0)}(\exists t>1:\hat\rho^{(0)}_t=1/2|\hat\rho^{(0)}_1=2,
\exists 0<s<1:\hat\rho^{(0)}_s=1/2)\\\nonumber
&&=
P_\rho^{(0)}(\exists t>1:\rho^{(0)}_t=1/2|\rho^{(0)}_1=2)<1
\end{eqnarray}
by transience of $\rho^{(0)}_.$, while
$$
\hat P_\rho^{(0)}(\exists t>1:\rho^{(0)}_t=1/2|\rho^{(0)}_1=2,
\nexists 0<s<1:\rho^{(0)}_s=1/2)=1
$$
since $T^{(0)}<\infty$ almost surely for $\hat \rho^{(0)}_.$.
  However, it is
immediately checked that the {\sl stopped} process which equals $\hat
\rho^{(s)}_t$ for $s\le t< T^{(s)}$ and, say, $0$ for $t\ge T^{(s)}$
is again Markovian.  This will play a role later.

We choose the parameter of the Bessel process as
$\delta=2(1+\alpha+\epsilon)$, with $\epsilon>0$ (this is the same
$\epsilon$ which appears in the statement of Theorem \ref{th:main}).
Then, from Eqs.  \eqref{eq:Kbess1}, \eqref{eq:Kbessel} and
\eqref{eq:K} it is immediate to realize that there exists
$p=p(\epsilon)$ with $0<p<1$ such that, for every $n\in \N$,
\begin{eqnarray}
  \label{eq:decomposK}
  K(n)=p\mathcal K^{(2(1+\alpha+\epsilon))}(n)+(1-p)\hat K(n)
\end{eqnarray}
where $\hat K(n)\ge 0$ and, of course, $\sum_{n\in\N}\hat
K(n)=1$. The important point here is the non-negativity of $\hat K(n)$, which
implies that both $\mathcal K(.)$ and
$\hat K(.)$ are probabilities on $\N$,  to which renewal processes are
naturally associated.

Note for later convenience that, as a consequence of \eqref{eq:Kstima},
\begin{eqnarray}
  \label{eq:lowBk}
  \frac{\mathcal K^{(2(1+\alpha+\epsilon))}(n)}{K(n)}\ge
  \frac{d_3(\epsilon)}{n^{2\epsilon}}.
\end{eqnarray}

\begin{rem}\rm
\label{rem:X}
Note that, if the slowly varying function $L(n)$ in \eqref{eq:K} tends
to a positive constant for $n\to\infty$, one can choose $\epsilon=0$
and in that case
\eqref{eq:lowBk} can be improved into
\begin{eqnarray}
\label{eq:lowBkimp}
\inf_{n\in\N} \frac{\mathcal K^{(2(1+\alpha))}(n)}{K(n)}\ge {d'_3}>0.
\end{eqnarray}
\end{rem}

\begin{figure}[h]
\begin{center}~
\leavevmode
\epsfxsize =14 cm
\psfragscanon
\psfrag{p}[c]{{ $\phi^{(x)}_t$}}
\psfrag{0}[c]{{ $0$}}
\psfrag{1}[c]{{ $1$}}
\psfrag{t}[c]{{ $t$}}
\psfrag{s}[c]{{ $\psi^{(x)}_t$}}
\psfrag{u}[c]{{ $1/2$}}
\psfrag{t1}[c]{{ $\tau_1$}}
\psfrag{t2}[c]{{ $\tau_2$}}
\psfrag{t3}[c]{{ $\tau_3$}}
\psfrag{yy}[c]{{ $T^{(0)}$}}
\psfrag{ii}[c]{{ $T^{(\tau_2)}$}}
\epsfbox{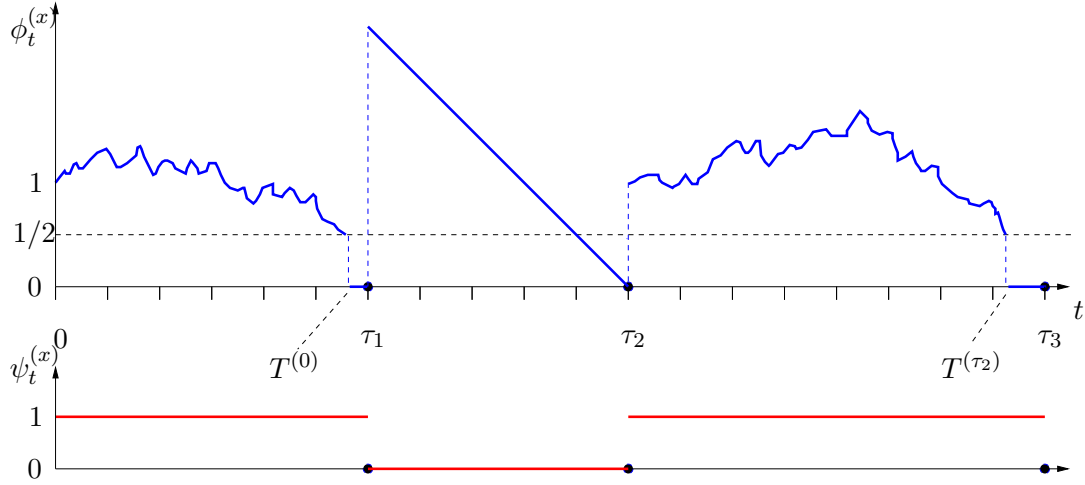}
\end{center}
\caption{\label{fig:markov}
An example of trajectory of $S^{(x)}_t=(\phi^{(x)}_t,\psi^{(x)}_t)$.
  In this picture the starting time $x$ equals $0$.  The top curve
  represents $\phi^{(x)}_t$, the bottom one $\psi^{(x)}_t$. In this
  example, $\psi^{(x)}_t=1$ between $0$ and $\tau_1$.  At the same
  time, $\phi^{(x)}_t$ performs a Bessel excursion starting from the
  value $1$, up to the time $T^{(0)}$ when it reaches the value $1/2$.
  Then it equals $0$ up to $\tau_1= \tilde T^{(0)}$. In the time
  interval $[\tau_1,\tau_2]$, on the other hand, $\psi^{(x)}_t$ equals
  $0$ and $\phi^{(x)}_t$ decreases linearly.  In the third time
  interval, one has again a Bessel excursion for $\phi^{(x)}$ and the value
  $1$ for $\psi^{(x)}$, and so on. The stretches of the trajectory
  $(\phi^{(x)}_t,\psi^{(x)}_t)$ between $\tau_i$ and $\tau_{i+1}$ are
  independent.  }
\end{figure}

Now, given $x\in\Z$ we construct a continuous-time Markov process
$\{S^{(x)}_t\}_{t\ge x}=\{(\phi^{(x)}_t,\psi^{(x)}_t)\}_{t\ge x}$,
with $\phi^{(x)}_t\ge 0$, $\psi^{(x)}_t\in \{0,1\}$ and initial
condition $S^{(x)}_x=(0,0)$. The process will satisfy the following
two properties:
\begin{itemize}
\item Let $t\in\Z$. Conditionally on $\phi^{(x)}_t=0$,
$\{S_u\}_{u>t}$ is independent of $\{S_u\}_{u<t}$.
\item Let $t_1<t_2\in \Z$. The process $\{S_u\}_{u>t_1}$, conditioned
on $\phi^{(x)}_{t_1}=0$, has the same law as $\{S_u\}_{u>t_2}$ conditioned
on $\phi^{(x)}_{t_2}=0$ and time-shifted to the left of $t_2-t_1$.
\end{itemize}
Therefore, we need to construct the trajectories only between two
successive integer times where $\phi^{(x)}_t=0$. The construction
proceeds as follows: whenever the condition
\begin{eqnarray}
  \label{eq:condizione}
  t\in \Z, \phi_t^{(x)}=0
\end{eqnarray}
is realized, we extract (independently of $\{S^{(x)}_u\}_{u\le t}$) a
random variable $\Psi$ which takes value $0$ with probability $(1-p)$,
and $1$ with probability $p$ ($p$ being defined in
Eq. \eqref{eq:decomposK}).  At that point (see Figure
\ref{fig:markov}):
\begin{itemize}
\item If $\Psi=0$, then we
extract a random variable $m\in\N$ with probability
law $\hat K(.)$ and we let
   $\phi^{(x)}_u=m+t-u$ for $u\in
  (t,t+m]$. In the same time interval, we let $\psi^{(x)}_u=\Psi=0$. At
  time $t+m$, we are back to condition \eqref{eq:condizione} and we
  start again the procedure with an independent extraction of $\Psi$.
\item If $\Psi=1$, then we let $\phi^{(x)}_u$ evolve like the process
  $\hat\rho^{(t)}_u$ for $u\in (t,t+T^{(t)})$ where, we recall,
  $T^{(t)}$ is the (random, but almost surely finite) first time after
  $t$ when $\hat\rho^{(t)}$ equals $1/2$. In particular,
  $\phi^{(x)}_{t^+}=1$.  Let $\tilde T^{(t)}=\inf\{j\in \Z:j\ge
  T^{(t)}\}$. Then, we let $\phi^{(x)}_u=0$ for $u\in [T^{(t)},\tilde
  T^{(t)}]$ and $\psi^{(x)}_u=\Psi=1$ for $u\in (t,\tilde T^{(t)}]$.  At
  time $\tilde T^{(t)}$ we are back to condition \eqref{eq:condizione}
  and we start again with an independent extraction of $\Psi$.
\end{itemize}
The process $S^{(x)}_.$ so constructed (whose law will be denoted by
$\hat \bP_x$), satisfies the following properties which are easily
checked:
\begin{enumerate}
\item If $\tau^{(x)}:=\{\Z\ni t\ge x: \phi^{(x)}_t=0\}$, then the
marginal distribution of $\tau^{(x)}$ is the law $\bP_x$ of Section
\ref{sec:notations} (the original renewal process associated to $K(.)$
with $\tau_0=x$). This is obvious from \eqref{eq:decomposK} and from
the construction of $S^{(x)}_.$.
\item Let
\begin{eqnarray}
  \label{eq:Boltz2}
  \frac{\dd\hat \bP_{x,y,\go}}{\dd\hat \bP_{x}}(S^{(x)}_.)=
\frac{e^{\sum_{n=x+1}^y(\beta \go_n-h)\ind_{\{n\in\tau^{(x)}\}}}
}
{\hat Z_{x,y,\go}}\ind_{\{y\in \tau^{(x)}\}}.
\end{eqnarray}
Then, the marginal distribution of $\tau^{(x)}$ is
the law $\bP_{x,y,\go}$ introduced in Eq. \eqref{eq:Boltz}.
\item For $(\beta,h)\in\cL$, the limit $\hat\bP_{\infty,\go}(f)$
obtained as $x\to- \infty,y\to\infty$ exists for every bounded local
observable $f$ (i.e., bounded function of $\{S^{(x)}_u\}_{u\in I}$,
$I$ bounded subset of $\R$.) This is a consequence of the fact that in
the localized region $\tau$ has a non-zero density in $\Z$ and that
the limit exists for functions depending only on $\tau$, as discussed
in Section \ref{sec:notations}.  We will call simply
$S_.=(\phi_.,\psi_.)$ the limit process obtained as $x\to-
\infty,y\to\infty$, and $\tau=\{t\in \Z:\phi_t=0\}$.
\item The process $S_.$ is Markovian. More precisely: if $A$ is a
local event supported on $[u,\infty)$ then
\begin{eqnarray}
  \label{eq:MarkovS}
  \hat \bP_{\infty,\go}(A|\{S_t\}_{t\le u})=\hat \bP_{\infty,\go}(A|S_u).
\end{eqnarray}
(This property is easily checked for $x,y$ finite, and then passes
to the thermodynamic limit).
\item
Let again $\tau=\{t\in \Z:\phi_t=0\}$ and $A_{a,b}$ the event
$\{a\in\tau, b\in \tau, \{a+1,\ldots,b-1\}\cap\tau=\emptyset\}$, for
$a,b\in\Z$ with $x<a<b<y$.  Under the law $\hat \bP_{x,y,\go}$,
conditionally on $A_{a,b}$, the variable $\psi_{a+}(=\psi_u$ for every
$u\in (a,b]$, from our construction of $S_.$) is independent of
$\{S_t\}_{t\in (-\infty,a)\cup (b,\infty)}$ and is a Bernoulli
variable which equals $0$ with probability $$ (1-p)\frac{\hat
K(b-a)}{K(b-a)} $$ and $1$ with probability $$ p\frac{\mathcal
K^{(2(1+\alpha+\epsilon))}(b-a)}{K(b-a)}\ge
\frac{d_4(\epsilon)}{(b-a)^{2\epsilon}},
$$ where the lower bound follows from \eqref{eq:lowBk}.  As for
$\{\phi_u\}_{u\in (a,b]}$, conditionally on $A_{a,b}$ it is also
independent of $\{S_t\}_{t\in (-\infty,a)\cup (b,\infty)}$.  If in
addition we condition on $\psi_{a+}=0$, then $\phi_u=b-u$, while if we
condition on $\psi_{a+}=1$ then $\{\phi_u\}_{u\in (a,b]}$ has the same
law as a trajectory of $\rho^{(a)}_u$ {\sl conditioned on} $T^{(a)}\in
(b-1,b]$ up to (and excluding) time $T^{(a)}$, and $\phi_u=0$ in
$[T^{(a)},b]$.  This property extends to the limit
$x\to-\infty,y\to\infty$.

\end{enumerate}

\section{The coupling inequality}

\label{sec:couplineq}

Consider two independent copies $S^1_.,S^2_.$ of the process $S_.$,
distributed according to the product measure $\hat \bP^{\otimes
  2}_{\infty,\go}(.)$.  As a consequence of property C of Section
\ref{sec:markov}, we can rewrite
\begin{eqnarray}
\label{eq:riscrittura}
  \bP_{\infty,\go}(k\in \tau|0\in\tau)-\bP_{\infty,\go}(k\in \tau)=
\hat \bE^{\otimes 2}_{\infty,\go}\left(\left.
\ind_{\{\phi^{1}_{k}=0\}}-\ind_{\{\phi^{2}_{k}=0\}}\right|\phi^1_0=0
\right).
\end{eqnarray}
Given two trajectories of $S_.$, define their first {\sl coupling time after
time zero} as
\begin{eqnarray}
  \label{eq:coupling}
  \mathcal T(S^1,S^2):=\inf\{t\ge 0: S^1_t=S^2_t\}.
\end{eqnarray}
It is important to remark that we are not requiring $ \mathcal
T(S^1,S^2)$ to be an integer.  Then, from the Markov property of $S$
it is clear that the r.h.s. of \eqref{eq:riscrittura} equals
\begin{eqnarray}
\hat \bE^{\otimes 2}_{\infty,\go}\left(\left.\left(
\ind_{\{\phi^{1}_{k}=0\}}-\ind_{\{\phi^{2}_{k}=0\}}\right)
\ind_{\{\mathcal T(S^1,S^2)>k\}}
\right|\phi^1_0=0
\right).
\end{eqnarray}
Therefore, we conclude that
\begin{eqnarray}
\label{eq:c}
  \left|\bP_{\infty,\go}(k\in \tau|0\in\tau)-\bP_{\infty,\go}(k\in
\tau)\right| \le \hat \bP^{\otimes 2}_{\infty,\go}\left(\left.\mathcal
T(S^1,S^2)>k\right|\phi^1_0=0 \right).
\end{eqnarray}
To proceed with the proof of Theorems \ref{th:main} and
\ref{th:mainAS} we are left with the task of giving upper bounds for
the probability that the coupling time is large. This will be done in
Section \ref{sec:largeT}, but first we need results on the geometry of
the set $\{t\in \Z:\phi_t=0\}\cap\{1,\ldots,k\}$, for $k$ large and
close to the critical line.

\begin{rem}\rm
\label{rem:decadAB}
In analogy with Eqs. \eqref{eq:riscrittura}-\eqref{eq:c}, under the
assumptions of Theorem \ref{th:AB} on the local observables $A,B$, one
has
\begin{eqnarray}
 \nonumber
\left|\bE_{\infty,\go}(A B)-\bE_{\infty,\go}(A)\bE_{\infty,\go}( B)\right|
&=&\left|\hat\bE^{\otimes 2}_{\infty,\go}\left[(A(\tau^1)B(\tau^1)-
A(\tau^1)B(\tau^2))
\ind_{\{\mathcal T(S^1,S^2)\ge k\}}\right]\right|\\
&\le&2||A||_\infty ||B||_\infty\hat \bP^{\otimes
2}_{\infty,\go}\left(\mathcal T(S^1,S^2)\ge k\right).
\end{eqnarray}
The upper bounds of Section \ref{sec:largeT} on the probability of
large coupling times imply therefore Theorem \ref{th:AB} (indeed, the
proof of Eqs. \eqref{eq:s1} and \eqref{eq:Tas} can be easily repeated
in absence of the conditioning on the event $\phi^1_0=0$.)
\end{rem}

\section{Estimates on the distribution of returns in a long time interval}

\label{sec:ritorni}

Ideas similar to those employed in this section have been already used
in Ref. \cite{cf:GT_ALEA} and, more recently, in \cite{cf:A}.

To simplify notations, we will from now on set $\bv:=(\beta,h)$,
$\mu:=\mu(\bv)$ and $\tf:=\tf(\bv)$. Also, in the following whenever a
constant $c(\bv)$ is such that for every bounded $B\subset \cL$ one
has $0<c_-(B)\le \inf_{\bv\in B}c(\bv) \le \sup_{\bv\in B}c(\bv)\le
c_+(B)<\infty$, we will say with some abuse of language that it is
independent of $\bv$.  In particular, this means that $c(\bv)$ cannot
vanish or diverge when the critical line is approached.

In this section we prove, roughly speaking, that if the interval
$\{1,\ldots,k\}$ is large there are sufficiently many points of $\tau$
in it, and that these points are rather uniformly distributed.  More
precisely: take the interval $\{1,\ldots,k\}$ and divide it into
disjoint blocks $B_\ell:=\{(\ell-1)R+1,\ldots,\ell R\}$,
$\ell=1,\ldots,M$ of size
\begin{eqnarray}
  \label{eq:R}
R:=\frac{c|\log \mu|}\mu,
\end{eqnarray}
where $c$ is a large (but independent of
$\bv$) positive constant to be chosen later and
\begin{eqnarray}
  \label{eq:M}
M=k\frac{\mu}{c|\log \mu|}.
\end{eqnarray}
In order to avoid a plethora of $\lfloor.\rfloor$, we are
assuming that $R$ and $M$ are integers. Let $\eta$ be a positive constant,
which will be chosen small (independently of $\bv$) later.
 Now we want to say that, with
probability at least $\simeq (1-\exp(-\mu k))$, a finite fraction of
the blocks contain at least a point of $\tau$:
\begin{proposition}
\label{prop:aver}
There exists $c_5<\infty$ such that
\begin{eqnarray}
\label{eq:prop61}
\bbE\, \bP_{\infty,\go}(\exists I\subset \{1,\ldots,M\}: |I|\ge \eta
M\;\mbox{\rm and}\; B_\ell\cap \tau=\emptyset \;\mbox{\rm for every}\;
\ell\in I)\le c_5 \mu^{-c_5}e^{-k \eta\,\mu/c_5}.
\end{eqnarray}
\end{proposition}

We will need also an analogous $\bbP(\dd\go)$-almost sure
result. However, in this case the strategy has to be modified and
$\{1,\ldots,k\}$ has to be divided into blocks whose lengths depend on
$\go$: namely, let $i_0(\go)=0$,
$$ i_j(\go)=\inf\{r> i_{j-1}(\go): Z_{i_{j-1}(\go),i_j(\go),\go}\ge
\frac1{\tf^c}\}
$$
and $M(\go)=\sup\{j:i_j(\go)\le k\}$.  Again, we define blocks $
B^\go_\ell:=\{i_{\ell-1}(\go)+1,\ldots,i_{\ell}(\go)\},
\ell=1,\ldots,M(\go)$, while $
B^\go_{M(\go)+1}:=\{i_{M(\go)}(\go)+1,\ldots,k\}$.  Then, one has:
\begin{proposition}
\label{prop:typ}
There exists $k_0(\go,\bv)$, $\bbP(\dd\go)$-almost surely finite and
$c_6(\bv)>0$ such that for every $k\ge k_0(\go,\bv)$
  \begin{enumerate}
\item
    \begin{eqnarray}
      \label{eq:Mgo}
 M(\go)\ge k \frac{\tf}{2c|\log \tf|}.
    \end{eqnarray}
\item
\begin{eqnarray}
\nonumber
&&  \bP_{\infty,\go}\left(\exists I\subset \{1,\ldots,M(\go)+1\}:
|I|\ge \eta M(\go)  \;\mbox{\rm and}\;
 B^\go_\ell\cap \tau=\emptyset\;\mbox{\rm for every}\;
\ell\in I\right)\\
\label{eq:ritornityp}
&&\le c_6(\bv)e^{-k \eta\,\tf/8}.
\end{eqnarray}
  \end{enumerate}
\end{proposition}

{\sl Proof of Proposition \ref{prop:aver}}
Define the event
$$ A:=\{\exists I\subset \{1,\ldots,M\}: |I|\ge \eta M\;\mbox{\rm
and}\; B_\ell\cap \tau=\emptyset \;\mbox{\rm for every}\; \ell\in I\}.
$$
Write
\begin{eqnarray}
  \label{eq:decomposiz}
&&  \bbE\, \bP_{\infty,\go}(A)=
\sumtwo{I\subset \{1,\ldots,M\}:}{|I|\ge \eta M}
\bbE \,\bP_{\infty,\go}(A_I)
\end{eqnarray}
where $A_I$ is the event
\begin{eqnarray}
\label{eq:AI}
A_I:=\{
B_\ell\cap \tau=\emptyset \;\mbox{\rm for every}\;
\ell\in I\}\cap\{B_\ell\cap \tau\ne\emptyset \;\mbox{\rm for every}\;
\ell\notin I\}
\end{eqnarray}
We can
rewrite (in a unique way) $B_I:=\cup_{\ell\in I} B_\ell$ as a disjoint
union of intervals,
\begin{eqnarray}
  \label{eq:disjoint}
B_I=\cup_{r=1}^{m(I)}\{i_r,\ldots,j_r\},
\end{eqnarray}
with $i_r\ge j_{r-1}+R$. In other words, any two adjacent blocks
$B_\ell,B_{\ell+1}$ with $\ell,\ell+1$ belonging to $I$ will be
regrouped in the same interval. Of course, $1\le m(I)\le |I|$ if $I$
is not empty.  Conditioning on the location $x_r$ of the first point
of $\tau$ at the left of $i_r$ and on the location $y_r$ of the first
point of $\tau$ at the right of $j_r$ one has
\begin{eqnarray}
\label{eq:spiegaz}
\bP_{\infty,\go}(A_I)\le e^{m(I)(|h|+\beta \go_{max})}\sumtwo{x_1\le i_1}
{j_1\le y_1\le (j_1+R)}
\sumtwo{(i_{m(I)}-R)\le x_{m(I)}\le i_{m(I)}}{y_{m(I)}\ge j_{m(I)}}
\sumthree{(i_r-R)\le x_r\le i_r}{j_r\le y_r\le (j_r+R)}{1< r<m(I)}
\prod_{r=1}^{m(I)}\frac 1{Z_{x_r,y_r,\go}}.
\end{eqnarray}
(If $m(I)=1$, the formula is slightly modified in that the sum is only
on $x_1\le i_1$ and $y_1\ge j_1$; the estimates which follow hold also
in this case).  Here we are using the fact that the disorder variables
are bounded, say, $|\go_n|\le \go_{max}$.  To obtain
\eqref{eq:spiegaz} observe that, if $i^-_r:=\max\{\tau_i:\tau_i\le
i_r\}$ and $j^+_r:=\min\{\tau_i:\tau_i\ge j_r\}$,
\begin{eqnarray}
 &&
\bP_{\infty,\go}(A_I;i_r^-=x_r,j^+_r=y_r\,\forall r=1,\ldots,m(I))\\
&&\le
\bP_{\infty,\go}(A_I|i_r^-=x_r,j^+_r=y_r\,\forall r=1,\ldots,m(I))
\le \prod_{r=1}^{m(I)}
\frac{K(y_r-x_r)e^{\beta \omega_{y_r}-h}}{Z_{x_r,y_r,\go}}
\end{eqnarray}
where we used \eqref{eq:identita} in the last step.  It is clear that,
on the event $A_I$, $i^-_r\ge i_r-R $ if $r>1$ (otherwise the block
$\{i_r-R,\ldots,i_r-1\}$ would be contained in $B_I$, which is not possible
due to $i_r\ge j_r+R$) and similarly
$j^+_r\le j_r+R$ if $r<m(I)$.  Then, \eqref{eq:spiegaz} immediately
follows.  Note that by the first inequality in \eqref{eq:ZZZ} one can
bound $Z_{x_r,y_r,\go}\ge
Z_{x_r,i_r,\go}Z_{i_r,j_r,\go}Z_{j_r,y_r,\go}$.  Therefore, using Eqs.
\eqref{eq:stimaZ}, \eqref{eq:Kstima} and
\eqref{eq:stimMu}, we get that
\begin{eqnarray}
\label{eq:conditioning}
\bbE\, \bP_{\infty,\go}(A_I)\le \mu^{-c_7}
\prod_{r=1}^{m(I)}c_7 R^{c_7} \bbE \frac 1{Z_{i_r,j_r,\go}}\le\mu^{-c_7}
\prod_{r=1}^{m(I)}c_7 R^{c_{7}}e^{-\mu(j_r-i_r)}(j_r-i_r)^{c_8}
\end{eqnarray}
for some positive $c_7,c_{8}$.
The factor $\mu^{-c_7}$ comes, through \eqref{eq:stimMu}, from the sum
$$ \sum_{x_1:x_1\le i_1}\bbE
\frac1{Z_{x_1,i_1,\go}}\left(=\sum_{y_{m(I)}:y_{m(I)}\ge j_{m(I)}}\bbE
\frac1{Z_{j_{m(I)},y_{m(I)},\go}}\right).
$$
Since $m(I)\le |I|$, one finds then
\begin{eqnarray}
  \label{eq:c10}
  \bbE\, \bP_{\infty,\go}(A_I)\le \mu^{-c_7}e^{-|I|(\mu R
-c_{7}\log R-\log c_7)}e^{c_8 \sum_{r=1}^{m(I)}\log (j_r-i_r)}.
\end{eqnarray}
Now we use Jensen's inequality for the logarithm and
the monotonicity of
 $x\rightarrow x\log (1/x)$ for $x>0$ small to bound
$$
e^{c_8 \sum_{r=1}^{m(I)}\log (j_r-i_r)} \le e^{c_8 |I|\log\left(
\frac k{|I|}\right)}.
$$
 From the definition of $R$ one sees then that, for $c$ sufficiently
large (independently of $\bv$)
 \begin{eqnarray}
  \bbE\,  \bP_{\infty,\go}(A_I)\le c_9\mu^{-c_7}
\exp\left(-\frac{c|I||\log \mu|}{2}\right)e^{c_8 |I|\log\left(
\frac k{|I|}\right)}
\end{eqnarray}
uniformly in $I$. Finally we can go back to the decomposition
\eqref{eq:decomposiz} which, together with elementary combinatorial
considerations, gives
\begin{eqnarray}
  \label{eq:finalAVE}
  \bbE\, \bP_{\infty,\go}(A)&\le& c_9\mu^{-c_7}\sum_{j\ge \eta M}
\left(
  \begin{array}{c}
M\\
j
  \end{array}
\right)e^{-j c |\log \mu|/2}e^{c_8 j \log \left(\frac{c\log \mu}{\eta\mu}
\right)}\\\nonumber
&\le&
c_{10}\mu^{-c_7}\left(
\begin{array}{c}
  M\\
M/2
\end{array}
\right)e^{-\eta k \mu/4}
\le
c_{11}\mu^{-c_7}e^{-\frac{\eta k \mu}8}
\end{eqnarray}
if  $c$ is large enough.
\hfill $\stackrel{{\small \mbox{Proposition \ref{prop:aver}}}}{\Box}$

{\sl Proof of Proposition \ref{prop:typ}} Observe first of all that,
thanks to \eqref{eq:ZZZ} and to the boundedness of disorder, for every
$\go$ and $x<y$
\begin{eqnarray}
  \frac1{c_{12}}\le \frac{Z_{x,y,\go}}{Z_{x,y+1,\go}}\le c_{12}
\end{eqnarray}
so that, say,
\begin{eqnarray}
\frac 1{\tf^c}\le   Z_{i_j(\go),i_{j+1}(\go),\go}\le \frac c{\tf^c}
\end{eqnarray}
if $c$ is sufficiently large (the lower bound holds by definition of
$i_j(\go)$, while the upper bound simply says that, since by definition
$Z_{i_j(\go),i_{j+1}(\go)-1,\go}<\tf^{-c}$, then
$Z_{i_j(\go),i_{j+1}(\go),\go}$ cannot be much larger than $\tf^{-c}$).
Therefore, denoting (with some abuse of notation) $i_{M(\go)+1}=k$ and
using repeatedly Eq. \eqref{eq:ZZZ}, we find
\begin{eqnarray}
  Z_{0,k,\go}\le \left(\frac c{\tf^c}\right)^{M(\go)+1}c_1^{M(\go)}
\prod_{r=1}^{M(\go)+1}(i_j(\go)-i_{j-1}(\go))^{c_1}
\end{eqnarray}
and, applying Jensen's inequality to the concave function
$x\rightarrow \log x$,
\begin{eqnarray}
\label{eq:pizza}
  \frac1k \log Z_{0,k,\go}\le c\frac{M(\go)+1}k |\log \tf|+(\log
c_1+\log c)\frac{M(\go)}k +c_1\frac{M(\go)+1}k \log \left(\frac k
{M(\go)+1} \right).
\end{eqnarray}
Now assume that
\begin{eqnarray}
  \label{eq:assurdo}
  \frac{M(\go)+1}k\le \frac{\tf}{2c|\log \tf|}.
\end{eqnarray}
Since the function $x\rightarrow x\log (1/x)$
is increasing for $x>0$ small, one deduces from \eqref{eq:pizza}
\begin{eqnarray}
  \frac1k \log Z_{0,k,\go}\le \frac34 {\tf}
\end{eqnarray}
if $c$ is chosen sufficiently large. But we know that $(1/k) \log
Z_{0,k,\go}$ converges to $\tf$ almost surely, and therefore
 the event \eqref{eq:assurdo} does not happen for $k$
larger than some random but finite  $k_0(\go)$. Equation \eqref{eq:Mgo}
is then proven.

As for \eqref{eq:ritornityp}, in view of Lemma \ref{corolla}
it is sufficient to prove that
\begin{eqnarray}
  \bP_{\infty,\go}\left(A;\{0,k+1\}\subset\tau\right)
\le c_6(\bv)e^{-k \eta\,\tf/8}
\end{eqnarray}
for $k\ge k_0(\go)$, where
$$
A^\go=\{\exists I\subset \{1,\ldots,M(\go)+1\}:
|I|\ge \eta M(\go)  \;\mbox{\rm and}\;
 B^\go_\ell\cap \tau=\emptyset\;\mbox{\rm for every}\;
\ell\in I\}.
$$
In analogy with Eqs. \eqref{eq:AI},
\eqref{eq:disjoint} define for $I\subset\{1,\ldots,M(\go)+1\}$
\begin{eqnarray}
 A^\go_I:=\{
B^\go_\ell\cap \tau=\emptyset \;\mbox{\rm for every}\;
\ell\in I\}\cap\{B^\go_\ell\cap \tau\ne\emptyset \;\mbox{\rm for every}\;
\ell\notin I\}
\end{eqnarray}
and rewrite $B_I:=\cup_{\ell\in I}B^\go_\ell$ as
$$
B_I=\cup_{r=1}^{m(I)}\{i_{x_r}(\go)+1,\ldots,i_{y_r}(\go)\}
$$
where the indices ${x_r},y_r$ are chosen so that $i_{x_r}(\go)\ge
i_{y_{r-1}}(\go)+2$.  Then, with a conditioning argument similar to the
one which led to Eq. \eqref{eq:conditioning}, one finds for $c$
sufficiently large
\begin{eqnarray}
\nonumber && \bP_{\infty,\go}( A^\go_I;\{0,k+1\}\subset\tau)\le
\bP_{\infty,\go}( A^\go_I|\{0,k+1\}\subset\tau)=
\bP_{0,k+1,\go}( A^\go_I)\\
&&\le
\tf^{c|I|}
\prod_{r=1}^{m(I)} c_{13}
[(i_{x_r}(\go)-i_{x_r-1}(\go))(j_{y_r+1}(\go)-j_{y_r}(\go))]^{c_{13}}\\
\nonumber && \le c_{14}^{|I|}e^{-c |I||\log \tf|}
\exp\left(c_{14}m(I)\log \left(\frac k{m(I)}\right)\right) \le
c_{15}(\bv)e^{-\frac c2|I||\log \tf|}.
\end{eqnarray}
In the third inequality we used, once more, Jensen's inequality for
 the logarithm function and in the fourth one the monotonicity of
 $x\rightarrow x\log (1/x)$ for $x>0$ small, plus Eq. \eqref{eq:Mgo}
 and the assumption that $|I|\ge \eta M(\go)$.  Considering all
 possible sets $I$ of cardinality not smaller than $\eta M(\go)$, we
 see that the l.h.s. of \eqref{eq:ritornityp} is bounded above by
\begin{eqnarray}
 c_{15}(\bv) \sum_{j\ge \eta M(\go)}\left(
  \begin{array}{c}
M(\go)+1\\
j
  \end{array}
\right)e^{-c j|\log \tf|/2}
\end{eqnarray}
and recalling \eqref{eq:Mgo}, the desired result Eq.
\eqref{eq:ritornityp} holds.
\hfill $\stackrel{{\small    \mbox{Proposition \ref{prop:typ}}}}{\Box}$

\section{Upper bounds on the probability of large coupling times}

\label{sec:largeT}

Finally, we can go back to the problem of estimating from above the
$\hat\bP^{\otimes 2}_{\infty,\go}$-probability that the coupling
time is larger than $k$, cf. Section \ref{sec:couplineq}.  This will
conclude the proof of Theorems \ref{th:main}, \ref{th:mainAS} and
\ref{th:AB}.

\subsection{The average case}
We wish first of all to prove that
\begin{eqnarray}
\label{eq:s1}
\bbE\,\hat \bP^{\otimes
2}_{\infty,\go}\left(\left.\mathcal T(S^1,S^2)>k+1\right|\phi^1_0=0
\right)\le \frac1{C_1(\epsilon)\mu^{1/C_1(\epsilon)}}e^{-k\,
C_1(\epsilon)\mu^{1+\epsilon}}.
\end{eqnarray}
To this purpose observe that, if $\tau^a=\{t\in \Z:\phi^a_t=0\}$,
$a=1,2$,
\begin{eqnarray}
\label{eq:purpose}
\nonumber &&\bbE \, \hat \bP^{\otimes 2}_{\infty,\go}\left(\left.  \exists
I\subset \{1,\ldots,M\}: |I|\ge \eta M,\,
B_\ell\cap\tau^1=\emptyset\;\mbox{or}\;
B_\ell\cap\tau^2=\emptyset\;\forall \ell\in
I\right|\phi^1_0=0\right)\\ &&=:\bbE \,\hat  \bP^{\otimes
2}_{\infty,\go}(\left.U\right|\phi^1_0=0)\le 2c_5\mu^{-c_5}e^{-k
\eta\mu/c_5}.
\end{eqnarray}
This would be an immediate consequence of Proposition \ref{prop:aver}
if the conditioning on $0\in\tau^1$ were absent. However, the proof of
Proposition \ref{prop:aver} can be repeated exactly in presence of
conditioning, i.e., when the measure $\bP_{\infty,\go}(.)$ is replaced
by $\bP_{0,\infty,\go}(.):=\lim_{y\to\infty}\bP_{0,y,\go}(.)$ in
Eq. \eqref{eq:prop61}.  Therefore,
\begin{eqnarray}
\label{eq:quasifinmu}
  \bbE\,\hat \bP^{\otimes 2}_{\infty,\go}\left(\left.\mathcal
T(S^1,S^2)>k+1\right|\phi^1_0=0 \right) &\le& 2c_5\mu^{-c_5}e^{-k
\eta\mu/c_5}\\\nonumber
&&+
\bbE\,\hat \bP^{\otimes 2}_{\infty,\go}\left(\left.\mathcal
T(S^1,S^2)>k+1\right|U^c,\phi^1_0=0 \right),
\end{eqnarray}
where $U^c$ is the complementary of the event $U$.  On the other hand,
provided that $\eta$ is chosen sufficiently small (but independent of
$\bv$) it is obvious that if the event $U^c$ occurs there exist at
least, say, $M/10$ integers $1<\ell_i< M$ such that
$\ell_i>\ell_{i-1}+2$ and $B_{r}\cap\tau^a\ne\emptyset$, for every
$a\in\{1,2\}$ and $r\in\{\ell_i-1,\ell_i,\ell_i+1\}$. The condition
$\ell_i>\ell_{i-1}+2$ simply guarantees that any two triplets of
blocks of the kind $\{B_{\ell_i-1},B_{\ell_i},B_{\ell_i+1}\}$ are
disjoint for different $i$, a condition we will need later in this section.
\begin{figure}[h]
\begin{center}~
\leavevmode
\epsfxsize =12 cm
\psfragscanon
\psfrag{x}[c]{{ $x\in \tau^1$}}
\psfrag{y}[c]{{ $y\in \tau^2$}}
\psfrag{j}[c]{{ $j$}}
\psfrag{k}[c]{{ $k$}}
\psfrag{d}[c]{{ $\phi^a_t$}}
\psfrag{1}[c]{{ $1$}}
\psfrag{12}[c]{{ $1/2$}}
\psfrag{m}[c]{{ $t_m$}}
\psfrag{t}[c]{{ $t$}}
\epsfbox{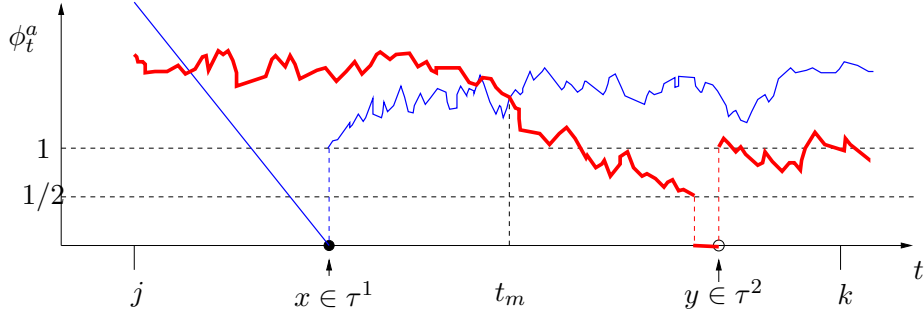}
\end{center}
\caption{\label{fig:good}
  An example of goodness. The thin line represents $\phi^1_t$ and the
  thick one represents $\phi^2_t$. The important thing is what happens
  between $x\in\tau^1$ and $y\in\tau^2$. Both paths perform a Bessel
  excursion in the time interval under consideration, which means that
  $\psi^1_t= \psi^2_t=1$ there. Since in this example $\phi^2_x>1$,
  there exists necessarily at least a time $t_m\in[x,y]$ where the
  two paths meet.  }
\end{figure}
We need to introduce the following definition:
\begin{definition}\rm
\label{def:goodblock}
A configuration of $(\tau^1,\tau^2)$ is called {\sl good in the interval
$\{j,\ldots,k\}$} if there exist $x,y\in \{j,\ldots,k\}$, with $x\le
y$, such that the following three conditions are satisfied:
\begin{itemize}
  \item either $\{x\in\tau^1\mbox{\;\;and\;\;}y\in \tau^2\}$ or
$\{x\in\tau^2\mbox{\;\;and\;\;}y\in \tau^1\}$
\item $\{x+1,\ldots,y-1\}\cap \tau^a=\emptyset$ for $a=1,2$
\item
$\psi^a_t=1$ for $a=1,2$ and $t\in [x,y]$.
\end{itemize}
\end{definition}
Roughly speaking (see Figure \ref{fig:good}), this means that
(assuming for definiteness $x\in\tau^1$) the point $x$ is overcome by
a Bessel excursion of $\phi^2_t$ which ends at $y$, while at $x$
starts a Bessel excursion of $\phi^1_t$ which overcomes $y$ and ends
at some later time.  Such a configuration is called {\sl good} in
$\{i,\ldots,j\}$ because the paths $S^1_t,S^2_t$ have a good chance of
meeting there, as the next result shows:
\begin{lemma}
\label{lemma:goodblock}
Conditionally on $(\tau^1,\tau^2)$ being good in the interval
$\{j,\ldots,k\}$ and on the configuration of $\{S^a_u\}_{u\notin [j,k]
}^{a=1,2}$, the $\hat \bP^{\otimes 2}_{\infty,\go}$-probability that
there exists $t\in [j,k]$ such that $S^1_t=S^2_t$ is bounded below by
a positive constant $c_0$, independent of $\go,j,k$ and of
$\{S^a_u\}_{u\notin [j,k]}^{a=1,2}$.
\end{lemma}
Lemma \ref{lemma:goodblock} is proven in Appendix \ref{sec:appBessel}.
Now recall property E, Section \ref{sec:markov}, of $\hat
\bP_{\infty,\go}$ and the discussion following Eq. \eqref{eq:purpose}
above, to conclude that, conditionally on the event $U^c$, the
configuration $(\tau^1,\tau^2)$ is good in each of the blocks
$B_{\ell_i}$ defined above, with probability at least $$
\left(\frac{d_4(\epsilon)}{R^{2\epsilon}}\right)^2.
$$ This holds {\sl independently} of what happens in $B_{\ell_j}$,
$j\ne i$, thanks to property \eqref{eq:identita}. Indeed, note that
there are points of $\tau^1$ and $\tau^2$ in both $B_{\ell_i-1}$ and
$B_{\ell_i+1}$ so that, via the Markov property, what happens in
$B_{\ell_i}$ is independent from what happens in $B_{\ell_j}$.

Using also Lemma \ref{lemma:goodblock} one has then that,
conditionally on $U^c$, the
$\hat \bP^{\otimes
2}_{\infty,\go}$-probability
that $\mathcal T(S^1,S^2)>k$ does not
exceed
\begin{eqnarray}
\label{eq:doesnot}
\left[1-c_0 \left(\frac{d_4(\epsilon)}{R^{2\epsilon}}\right)^2
\right]^{M/10}.
\end{eqnarray}
Recalling the definitions \eqref{eq:R} and \eqref{eq:M} of $R$ and $M$,
one can bound this probability from above with
\begin{eqnarray}
  \exp\left(-d_5(\epsilon)k \mu^{1+5\epsilon}/|\log \mu|^2\right).
\end{eqnarray}
Together with Eq. \eqref{eq:quasifinmu}, this concludes the
proof of Eq. \eqref{eq:risultaAV}.
\hfill $\stackrel{\mbox{\tiny Theorem \ref{th:main}, Eq.
\eqref{eq:risultaAV}}}{\Box}$

\begin{rem}\rm
\label{rem:XY}
A look at Remark \ref{rem:X} shows that, if the slowly varying
function $L(n)$ in \eqref{eq:K} tends to a positive constant for
$n\to\infty$, the upper bound \eqref{eq:doesnot} can be improved into
$ \exp(-\tilde d_5M)$ with $\tilde d_5>0$.  From this and
Eq. \eqref{eq:quasifinmu} the claim of Remark \ref{rem:extract}
follows immediately.
\end{rem}

\subsection{The almost-sure case}
Let us finally prove that, almost surely,
\begin{eqnarray}
\label{eq:Tas}
 \hat \bP^{\otimes
2}_{\infty,\go}\left(\left.\mathcal T(S^1,S^2)>k+1\right|\phi^1_0=0
\right)\le  C_2(\go)e^{-C_1\tf^{1+\epsilon}}.
\end{eqnarray}
The proof is quite similar to that of the average case. Define (with
the notations of Section \ref{sec:ritorni}) the event
\begin{eqnarray}
  \label{eq:W}
  W(\go):=\left\{\exists I\subset \{1,\ldots,M(\go)\}:|I|\ge \eta
k\frac{\tf}{2c|\log \tf|}, B_\ell^\go\cap
\tau^1=\emptyset\;\mbox{or}\; B^\go_\ell\cap\tau^2=\emptyset\;\forall
\ell\in I\right\}
\end{eqnarray}
so that
\begin{eqnarray}
 \label{eq:77}
 \hat \bP^{\otimes
2}_{\infty,\go}\left(\left.\mathcal T(S^1,S^2)>k+1\right|\phi^1_0=0
\right)&\le& 2c_6(\bv)e^{-k\eta\,\tf/16}\\\nonumber
&&+ \hat \bP^{\otimes
2}_{\infty,\go}\left(\left.\mathcal T(S^1,S^2)>k+1\right|W(\go)^c,\phi^1_0=0
\right)
\end{eqnarray}
$\bbP(\dd\go)$-almost surely, for $k>k_0(\go)$. If the event
$W(\go)^c$ occurs, one can find $G(\go)\ge k\tf/(20c|\log \tf|)$
integers $1<\ell_i< M$ such that $\ell_i>\ell_{i-1}+2$ and
$B^\go_{r}\cap\tau^a\ne\emptyset$, for every $a\in\{1,2\}$ and
$r\in\{\ell_i-1,\ell_i,\ell_i+1\}$. $(\tau^1,\tau^2)$ is good in each
of the blocks $B^\go_{\ell_j}$ with probability at least $$
\left(
\frac{d_4(\epsilon)}{(i_{\ell_j}(\go)-i_{\ell_j-1}(\go))^{2\epsilon}}
\right)^2.
$$
Therefore, conditionally on $W(\go)^c$, the
$\hat \bP^{\otimes
2}_{\infty,\go}$-probability that $\mathcal T(S^1,S^2)>k$ does not
exceed
\begin{eqnarray}
&& \prod_{j=1}^{G(\go)}\left[1-c_0\left(
\frac{d_4(\epsilon)}{(i_{\ell_j}(\go)-i_{\ell_j-1}(\go))^{2\epsilon}}
\right)^2\right]\le \exp\left( -d_6(\epsilon)\sum_{j=1}^{G(\go)}
(i_{\ell_j}(\go)-i_{\ell_j-1}(\go))^{-4\epsilon} \right)\\\nonumber
&&\le
\exp\left[-d_6(\epsilon)G(\go)\left(\frac{G(\go)}{\sum_{j=1}^{G(\go)}
(i_{\ell_j}(\go)-i_{\ell_j-1}(\go))}\right)^{4\epsilon}\right] \le
\exp\left[-d_6(\epsilon)k\left(\frac{G(\go)}{k}\right)^{1+4\epsilon}\right],
\end{eqnarray}
where we used Jensen's inequality for the convex function
$x\rightarrow x^{-4\epsilon}$.
The lower bound $G(\go)\ge k\tf/(20c|\log \tf|)$, together with
Eq. \eqref{eq:77} are then enough to obtain the desired estimate
\eqref{eq:Tas}.

\hfill $\stackrel{\mbox{\tiny Theorem \ref{th:main},
Eq. \eqref{eq:risultaAS}}}{\Box}$

\appendix

\section{Some technical facts on Bessel processes}

\label{sec:appBessel}

Take $0<b<a<\infty$ and consider a Bessel process $\rho^{(0)}_t,
t\ge0$ of dimension $\delta$ starting from $\rho^{(0)}_0=a$ at time
$0$. Let $T_{a,b}$ be the first hitting time of $b$, i.e.,
$T_{a,b}=\inf\{t\ge0:\rho^{(0)}_t=b\}$.  Then, it follows from
\cite[Theorem 3.1]{cf:kent} plus \cite[Theorem 2.5]{cf:mourad} that,
conditionally on $T_{a,b}<\infty$, the density of the probability
distribution of $T_{a,b}$ with respect to the Lebesgue measure on
$\R^+$ is proportional to
\begin{eqnarray}
  \label{eq:Mourad}
p(t):=\int_0^\infty B(z)e^{-tz/2}dz:=\int_0^\infty
\frac{J_\nu(b\sqrt z)Y_\nu(a\sqrt z)-J_\nu(a\sqrt z)Y_\nu(b\sqrt z)
}
{J^2_\nu(b\sqrt z)+Y^2_\nu(b\sqrt z)}
e^{-tz/2}dz
\end{eqnarray}
where $J_\nu(z)$ and $Y_\nu(z)$ are Bessel function of the first and
second kind, respectively \cite[Chapter 7.2.1]{cf:erdely}, and
$\nu=(\delta/2)-1$.  From \cite[Chap. 7.2.1, Eqs. (3)-(4)]{cf:erdely}
one deduces that $B(z)z^{-\nu}\to c_\nu(a,b)$ for $z\to 0^+$, where
$c_\nu(a,b)$ is a finite and positive constant whose precise value is
not needed for our purposes.  Therefore, the Abelian Theorem
\cite[Chapter 5, Corollary 1a]{cf:widder} gives
\begin{eqnarray}
  \label{eq:asympt}
p(t)t^{\nu+1}=p(t)t^{\delta/2}\stackrel{t\to+\infty}\longrightarrow
\frac{c_\nu(a,b)\Gamma(\nu+2)2^{\nu+1}}{\nu+1}.
\end{eqnarray}
  From Eq.  \eqref{eq:asympt}, the asymptotic behavior
\eqref{eq:Kbessel} immediately follows taking $a=1,b=1/2$ (of course,
any other values $0<b<a<\infty$ would be equally good).

\subsection{Proof of Lemma \ref{lemma:goodblock}} Let $x,y$
be any pair of sites which satisfies the conditions required by
Definition \ref{def:goodblock}.  Assume for definiteness that
$x\in\tau^1,y\in \tau^2$.  We assume also that $x<y$, otherwise the
lemma is trivial.  For technical reasons, it is also convenient to
treat apart the case $x=y-1$. In this case, the lemma follows
immediately from
\eqref{eq:ZZZ}. Indeed, from this is easily deduced in particular
that, conditionally on $y\in\tau^2$, the probability that also
$y-1\in\tau^2$ is greater than some positive constant, independent of
$\go$.

As for the more difficult case where $x<y-1$, it is clear that there
exists $x\le t\le y$ such that $\phi^1_t=\phi^2_t$ whenever
$\phi^2_x\ge1$ (we assume that $x\ne \tau^2$, otherwise the existence
of $t$ such that $\phi^1_t=\phi^2_t$ is trivial).  This follows (see
also Figure \ref{fig:good}) from the observation that
$\phi^1_{x^+}=1,\phi^1_y\ge1/2$ and that there exists $y-1<s\le y$ with
$\phi^2_s=1/2$, together with the fact that the trajectories of the
Bessel process are continuous almost surely.  Therefore, the Lemma
follows if we can prove that the probability that $\phi^2_x\ge1$ is
bounded below by a positive constant.  This is the content of
\eqref{eq:supf} below.

In order to state \eqref{eq:supf}, we need to introduce the Bessel
Bridge process of dimension $\delta$ \cite[Chapter XI.3]{cf:RY}. Given
$u\ge0$ and $a,v>0$, the Bessel Bridge is a continuous process
$\{X_t\}_{t\in[0,a]}$ (whose law is denoted by $P^{a,\delta}_{u,v}$)
which starts from $u$ at time $0$, ends at $v$ at time $a$ and such
that, given $0<s_1<\ldots <s_k<a$, the law of
$(X_{s_1},\ldots,X_{s_k})$ has density
\begin{eqnarray}
\label{eq:bbridge}
p^\delta_{s_1}(u,x_1)p^\delta_{s_2-s_1}(x_1,x_2)\ldots p^\delta_{a-s_k}
(x_k,v)/p^\delta_a(u,v).
\end{eqnarray}
Then, what we need is
\begin{eqnarray}
  \label{eq:supf}
  \inf_{u,v\ge 1/2} P^{2,\delta}_{u,v}(X_1\ge1|X_s>1/2\; \forall\;
s\in [0,2])>0.
\end{eqnarray}
Of course, $u,v$ correspond to the values
$\phi^{2}_{x-1},\phi^{2}_{x+1}$, respectively. It is immediate to
realize that \eqref{eq:supf} concludes the proof of Lemma
\ref{lemma:goodblock}.

Inequality \eqref{eq:supf} is easily proven: indeed, via FKG
inequalities \cite{cf:FKG} \cite{cf:laroche} one has (see details below)
\begin{eqnarray}
\label{eq:dafkg}
P^{2,\delta}_{u,v}(X_1\ge1|X_s>1/2\; \forall\;
s\in [0,2])
\ge P^{2,\delta}_{u,v}(X_1\ge1).
\end{eqnarray}
Using formula \eqref{eq:bbridge}, the r.h.s. of \eqref{eq:dafkg}
equals
\begin{eqnarray}
 \frac{  \int_1^\infty p_1^\delta(u,w)p_1^\delta(w,v)\dd w}
{p^\delta_2(u,v)}=\frac{e^{-(u^2+v^2)/4}}{I_\nu(uv/2)}\int_1^\infty
w\,e^{-w^2}I_\nu(uw)I_\nu(vw)\dd w.
\end{eqnarray}
Since $I_\nu(w)>0$ for $w>0$ and
\begin{eqnarray}
  \label{eq:asymptBI}
  \lim_{w\to\infty}e^{-w}\sqrt w I_\nu(w)\in (0,\infty)
\end{eqnarray}
(this can be extracted from \cite[
Chap. 7.13.1, Eq. (5); cf. Chap. 7.2.6 for the definition of the
Hankel symbol $(\nu,m)$]{cf:erdely}), one has
\begin{eqnarray}
  0<c^-_\nu:=\inf_{z\ge1/2} e^{-z}\sqrt z I_\nu(z)\le \sup_{z\ge1/2}
  e^{-z}\sqrt z I_\nu(z)=:c_\nu^+<\infty.
\end{eqnarray}
Therefore, the l.h.s. of \eqref{eq:dafkg} is bounded below by
\begin{eqnarray}
 \frac{(c_\nu^-)^2}{\sqrt 2 c^+_\nu}\int_1^\infty
e^{-\left(w-\frac{u+v}2\right)^2}\dd w,
\end{eqnarray}
which tends to a positive constant if $u\to+\infty$ or $v\to+\infty$
(or both), thus yielding Eq. \eqref{eq:supf}.

Finally, we show how \eqref{eq:dafkg} follows from the FKG
inequalities. Due to the continuity of the trajectories of the Bessel
Bridge, the probability in the l.h.s. of \eqref{eq:supf} equals
\begin{eqnarray}
\label{eq:quasiFKG}
\lim_{n\to\infty}
P^{2,\gamma}_{u,v}(X_1\ge1|X_{i/n}>1/2,\;
i=1,\ldots,2n-1).
\end{eqnarray}
Let $p(x_1,\ldots,x_{2n-1})$ be the probability density of
$(X_{1/n},\ldots,X_{(2n-1)/n})$.  Given $\underline
x^a:=(x^a_1,\ldots,x^a_{2n-1})$, $x^a_j>0$, $a=1,2$, define
$\underline x^1\vee \underline x^2:=((x^1_1\vee
x^2_1),\ldots,(x^1_{2n-1}\vee x^2_{2n-1}))$ and analogously
$\underline x^1\wedge \underline x^2$. Then, from the continuity and
Markov property of the Bessel Bridge process \cite[Chapter
XI.3]{cf:RY} it is clear that $p(\underline x^1\vee \underline x^2)p(
\underline x^1\wedge \underline x^2)\ge p(\underline x^1)p(\underline
x^2)$. This is just the FKG inequality, which implies in particular
that the probability in \eqref{eq:quasiFKG}, for any given $n$, is not
smaller than $P^{2,\delta}_{u,v}(X_1\ge1)$.

\hfill $\stackrel{\mbox{\tiny Lemma \ref{lemma:goodblock}}}{\Box}$

\section{Technical estimates on $Z_{x,y,\go}$ and $\bP_{\infty,\go}$}

\label{app:Z}

In this section we collect some technical estimates, which in very similar
form have been already used in the previous literature.
Let us notice at first that, for every $x<y$ and uniformly in $\go$,
\begin{eqnarray}
  \label{eq:stimaZ}
  Z_{x,y,\go}\ge e^{\beta \go_{y}-h}K(y-x).
\end{eqnarray}
Also, Eq. \eqref{eq:K} and the property of slow
variation imply that for every $\epsilon>0$ there exist positive
constants $d_1(\epsilon),d_2(\epsilon)$ such that, for every $n\in
\N$,
\begin{eqnarray}
  \label{eq:Kstima}
\frac{d_1(\epsilon)}{n ^{1+\alpha+\epsilon}}\le
K(n)\le\frac{d_2(\epsilon)}{n^{1+\alpha-\epsilon}}.
\end{eqnarray}

In Lemma A.1 of \cite{cf:GT_ALEA} it was proven that there exists
$c_1$, which in the case of bounded disorder can be chosen independent
of $\go$, such that for every $x<z<y$
\begin{eqnarray}
  \label{eq:ZZZ}
Z_{x,z,\go}Z_{z,y,\go}\le  Z_{x,y,\go}\le c_1 ((z-x)\wedge(y-z))^{c_1}
Z_{x,z,\go}Z_{z,y,\go}.
\end{eqnarray}
As it was shown in \cite[Proposition 2.7]{cf:GT_ALEA}, this
immediately implies that there exists $c'_1>0$ such that, for every
$y>x$,
$$
\left|\frac1{|y-x|}\bbE \log Z_{x,y,\go}-\tf(\bv)\right|\le
c'_1\frac{\log |y-x|}{|y-x|}.
$$
 Similarly, one can see that
\begin{eqnarray}
  \label{eq:stimMu}
 \left|-\frac1{|y-x|}\log \bbE\frac1{ Z_{x,y,\go}}-\mu(\bv)\right|\le
 c'_1 \frac{\log |y-x|}{|y-x|}:
\end{eqnarray}
this follows immediately observing that  \eqref{eq:ZZZ} implies
\begin{eqnarray}
\frac1{c_1(2N)^{c_1}} \left(\bbE \frac1{Z_{-N,N,\go}}\right)^2\le
\bbE \frac1{Z_{-2N,2N,\go}}\le \left(\bbE \frac1{Z_{-N,N,\go}}\right)^2.
\end{eqnarray}
A minor modification of the proof of \cite[Lemma A.1]{cf:GT_ALEA} gives
also
\begin{lemma}
\label{lemma1}
  Let $A$ be a local event supported in $\{1,\ldots,a\}$. Then,
  \begin{eqnarray}
    \label{eq:estimate}
    \bP_{\infty,\go}(A;\{-k,\ldots,0\}\cap \tau\ne\emptyset)\le c_1(a k)^{c_1}
\bP_{\infty,\go}(A;\{0,a+1\}\in \tau).
  \end{eqnarray}
\end{lemma}
We will also need the following result, which follows from \cite[Lemma
3.1]{cf:GT_ALEA}:
\begin{lemma}
\label{easylemma}
For every $\bv\in \mathcal L$ there exist positive constants
$c_2(\go;\bv),c_3 (\bv)$ (with $c_2(\go;\bv)$ finite
$\bbP(\dd\go)-$almost surely) such that, for every $k\in\N$,
\begin{eqnarray}
  \bP_{\infty,\go}(\tau\cap\{-k^2,\ldots,0\}=\emptyset)\le
c_2(\go;\bv)e^{-c_3(\bv) k^2}.
\end{eqnarray}
\end{lemma}
As a consequence of Lemmas \ref{easylemma} and
\ref{lemma1}, we have finally
\begin{lemma}
\label{corolla}
  Let $A$ be a local event supported in $\{1,\ldots,k\}$. Then,
\begin{eqnarray}
  \label{eq:cor1}
\bP_{\infty,\go}(A)\le c_1 k^{c_1}
\bP_{\infty,\go}(A;\{0,k+1\}\subset \tau) +c_2(\go;\bv)e^{-c_3(\bv)
k^2}.
\end{eqnarray}
\end{lemma}

\section*{acknowledgments}
This work originated from discussions with Giambattista Giacomin, to
whom I am very grateful for several suggestions. Partial
supported by the GIP-ANR project JC05\_42461 ({\sl POLINTBIO}) is
acknowledged.

\end{document}